\documentclass[onefignum,onetabnum]{siamart190516}

\usepackage{lipsum}
\usepackage{amsfonts}
\usepackage{graphicx}
\usepackage{epstopdf}
\usepackage{algorithm}
\usepackage{algpseudocode}
\usepackage{multicol}
\usepackage{siunitx}
\usepackage{subcaption}
\ifpdf
  \DeclareGraphicsExtensions{.eps,.pdf,.png,.jpg}
\else
  \DeclareGraphicsExtensions{.eps}
\fi
\usepackage{tikz}
\usetikzlibrary{calc}
\usetikzlibrary{quotes}
\usetikzlibrary{shapes}
\tikzstyle{vertexstyle}=[thick, draw, rectangle, fill opacity=0.4, text opacity=1]
\tikzstyle{edgestyle}=[thick, color=lightgray]

\usepackage{amsopn}
\DeclareMathOperator{\diag}{diag}
\DeclareMathOperator{\rank}{rank}

\usepackage{bm}
\newcommand{\vct}[1]{\bm{\mathsf{#1}}}
\newcommand{\mtx}[1]{\bm{\mathsf{#1}}}
\newcommand{\mA}{\mtx{A}}
\newcommand{\mB}{\mtx{B}}
\newcommand{\mD}{\mtx{D}}
\newcommand{\mG}{\mtx{G}}
\newcommand{\mI}{\mtx{I}}
\newcommand{\mP}{\mtx{P}}
\newcommand{\mQ}{\mtx{Q}}
\newcommand{\mR}{\mtx{R}}

\newcommand{\mU}{\mtx{U}}
\newcommand{\mV}{\mtx{V}}

\newcommand{\mY}{\mtx{Y}}
\newcommand{\mZ}{\mtx{Z}}
\newcommand{\mOmega}{\mtx{\Omega}}
\newcommand{\mPsi}{\mtx{\Psi}}
\newcommand{\mzero}{\mtx{0}}
\newcommand{\qq}{\vct{q}}
\newcommand{\uu}{\vct{u}}
\newcommand{\vv}{\vct{v}}
\newcommand{\xx}{\vct{x}}

\newcommand{\tree}{\mathcal{T}}                               %

\newcommand{\R}{\mathbb{R}}                                   %
\newcommand{\col}{\texttt{col}}                               %
\newcommand{\nullspace}{\texttt{null}}                        %

\newcommand{\bigO}{\mathcal{O}}

\newcommand{\tmult}{T_\text{mult}} %

\newcommand{\trand}{T_\text{rand}} %
\newcommand{\tflop}{T_\text{flop}}
\newcommand{\tcompress}{T_\text{compress}} %

\def\lrcolor{white!40!blue} %
\def\frcolor{white!30!red} %
\def\randcolor{lightgray} %

\newsiamremark{remark}{Remark}
\newsiamremark{hypothesis}{Hypothesis}
\crefname{hypothesis}{Hypothesis}{Hypotheses}
\newsiamthm{claim}{Claim}

\headers{
    Randomized Compression of Rank-Structured Matrices
}{J. Levitt and P. G. Martinsson}

\title{
    Linear-Complexity Black-Box Randomized Compression of Rank-Structured Matrices
    \thanks{Submitted to the editors DATE.
}}

\author{
    James Levitt\thanks{
        Oden Institute for Computational Engineering and Sciences,
        University of Texas at Austin, Austin, TX
        (\email{jlevitt@utexas.edu}, \email{pgm@oden.utexas.edu}).
    }
    \and
    Per-Gunnar Martinsson\footnotemark[2]
}

\newcommand{\pvct}[1]{#1}

\ifpdf
\hypersetup{
    pdftitle={
        Linear-Complexity Black-Box Randomized Compression of Rank-Structured Matrices
    },
    pdfauthor={J. Levitt and P. G. Martinsson}
}
\fi

\begin{document}

\maketitle

\begin{abstract}
A randomized algorithm for computing a compressed representation of a
given rank-structured matrix $\mtx{A} \in \mathbb{R}^{N\times N}$ is presented.
The algorithm interacts with $\mtx{A}$ only through its action on vectors.
Specifically, it draws two tall thin matrices
$\mtx{\Omega},\,\mtx{\Psi} \in \mathbb{R}^{N\times s}$ from a suitable distribution,
and then reconstructs $\mtx{A}$ from the information contained in the set $\{\mtx{A}\mtx{\Omega},\,\mtx{\Omega},\,\mtx{A}^{*}\mtx{\Psi},\,\mtx{\Psi}\}$.
For the specific case of a ``Hierarchically Block Separable (HBS)'' matrix
(a.k.a.~Hierarchically Semi-Separable matrix) of block rank $k$,
the number of samples $s$ required satisfies $s = O(k)$, with $s \approx 3k$
being representative.
While a number of randomized algorithms for compressing rank-structured
matrices have previously been published, the current algorithm appears
to be the first that is both of truly linear complexity (no $N\log(N)$ factors in the complexity bound) and fully ``black box'' in the sense that no matrix
entry evaluation is required.
Further, all samples can be extracted in parallel, enabling
the algorithm to work in a ``streaming'' or ``single view'' mode.
\end{abstract}



\section{Introduction} \label{sec:intro}

This work describes an efficient algorithm for handling
large dense matrices that have \textit{rank structure}. To simplify
slightly, this means that an $N\times N$ matrix can be tessellated
into $O(N)$ blocks in such a way that each block is either small
or of low numerical rank, cf.~\cref{fig:matrix_hodlr_lvl3}. This structure allows the matrix to be stored
and applied to vectors efficiently, often with cost that scales linearly
or close to linearly with $N$ \cite{hackbusch1999sparse,2010_borm_book,bebendorf2008hierarchical,2019_martinsson_fast_direct_solvers}. Sometimes, it is also possible to compute
an approximate inverse or LU factorization in linear or close to linear
time \cite{2016_vogel_xia_superfast_DC}. 
Matrices of this type have turned out to be ubiquitous in both
engineering and data sciences, and have been the subject of much research
in recent decades, going under names such as
$\mathcal{H}$-matrices~\cite{bebendorf2008hierarchical,2010_borm_book,hackbusch1999sparse};
HODLR matrices~\cite{ambikasaran2013mathcal,martinsson2009fast},
Hierarchically Block Separable (HBS) or Hierarchically Semi-Separable (HSS)
matrices~\cite{chandrasekaran2007fast,chandrasekaran2006fast,xia2010superfast},
Recursive Skeletonization~\cite{gillman2012direct,ho2012fast,2005_martinsson_fastdirect,minden2017recursive},
and many more.

The specific problem we address is the following: Suppose that $\mtx{A}$
is an $N\times N$ matrix that we know has HBS structure,
but we do not have direct access to the low-rank factors that define
the compressible off-diagonal blocks.
Instead, we have access to fast ``black-box'' algorithms
that given tall thin matrices
$\mtx{\Omega},\mtx{\Psi} \in \mathbb{R}^{N\times s}$,
evaluate the matrix-matrix products
$$
\mtx{Y} = \mtx{A\Omega}
\qquad\mbox{and}\qquad
\mtx{Z} = \mtx{A}^{*}\mtx{\Psi}.
$$
The problem is then to recover $\mA$
from the information in the set
$\{\mtx{Y},\,\mtx{\Omega},\,\mtx{Z},\,\mtx{\Psi}\}$.
The algorithm described here solves the reconstruction 
problem using $s = O(k)$ sample vectors, where $k$ is 
an upper bound on the ranks of the off-diagonal blocks
that we assume is known in advance.
(The scaling factor hidden in the formula $s = O(k)$ is modest,
with $s \approx 3k$ being representative.)

The scheme presented has several important applications. 
First, it can be
used to derive a rank-structured representation of any integral operator
for which a fast matrix-vector multiplication algorithm,
such as the Fast Multipole Method~\cite{greengard1987fast,greengard1997new},
is available. 
The new representation that we compute opens the door to a wider
range of matrix operations such as LU factorization, matrix inversion,
and sometimes even full spectral decompositions. 
Second, it can greatly
simplify algebraic operations involving products of rank-structured
matrices. For instance, the perhaps key application of rank-structured
matrix algebra is the acceleration of sparse direct solvers, as the
dense matrices that arise during LU factorization are often rank structured.
In the course of such a solver, a typical operation would
be to form a Schur complement such as
$\mtx{S}_{22} = \mtx{A}_{21}\mtx{A}_{11}^{-1}\mtx{A}_{12}$ that would
arise when the top left block $\mtx{A}_{11}$ is eliminated from a $2\times 2$ blocked
matrix. If $\mtx{A}_{11}$ is rank structured, then $\mtx{A}_{11}^{-1}$
can easily be applied to vectors via an LU factorization. If, additionally,
$\mtx{A}_{12}$ and $\mtx{A}_{21}$ are either sparse or rank structured,
then $\mtx{S}_{22}$ can easily be applied to a vector.
The technique described will then enable one to construct a data-sparse
representation of $\mtx{S}_{22}$.
In contrast, to directly evaluate the product $\mtx{A}_{21}\mtx{A}_{11}^{-1}\mtx{A}_{12}$
is both onerous to code and slow to execute.

The method we describe is inspired by the randomized compression algorithm
introduced in~\cite{2008_martinsson_randomhudson,martinsson2011fast} for
compressing HBS matrices.
That algorithm is similar to the one presented here in that it requires
only $O(k)$ matrix-vector applications and $\bigO(k^2 N)$ floating point operations.
However, the algorithm of \cite{2008_martinsson_randomhudson,martinsson2011fast}
is not a true black-box algorithm since, in addition to randomized samples,
it also requires direct evaluation of a small number of entries of the matrix.
In contrast, the method presented here is truly black box.

\begin{remark}[Peeling algorithms]
A related class of algorithms
for randomized compression of rank-structured matrices
is described in~\cite{lin2011fast,martinsson2016compressing},
and a recent improvement is described in~\cite{levitt2022randomized}.
These techniques are ``true'' black-box algorithms
in that they only access the matrix
through the black-box matrix-vector multiplication routines,
and they also apply to a broader class of rank-structured
matrices than the one considered here.
However, they require $O(k\,\log(N))$ samples
and $O(k^2\,N\,\log(N))$ floating point operations,
so they do not have linear complexity.
\end{remark}

The manuscript is structured as follows:
Section \ref{sec:preliminaries} surveys some basic linear algebraic techniques that
we rely on.
Section  \ref{sec:hss} introduces our formalism for
HBS matrices.
Section \ref{sec:compress} describes the new algorithm, and analyzes its
asymptotic complexity.
Section \ref{sec:experiments} describes numerical results.

\section{Preliminaries} \label{sec:preliminaries}
In introducing well-known material,
we follow the presentation of~\cite{martinsson2016compressing}.

\subsection{Notation}

Throughout the paper, we measure a vector $\xx \in \R^n$
by its Euclidean norm
$\Vert \xx \Vert = \left( \sum_i \vert x_i \vert^2 \right)^\frac{1}{2}$.
We measure a matrix $\mA \in \R^{m \times n}$
with the corresponding operator norm
$\Vert \mA \Vert = \sup_{\Vert \xx \Vert = 1} \Vert \mA \xx \Vert$,
and in some cases with the Frobenius norm
$\Vert \mA \Vert_\text{Fro} = (\sum_{i, j} \vert \mA(i, j) \vert^2)^{1/2}$.
To denote submatrices,
we use the notation of Golub and Van Loan~\cite{golub2013matrix}:
If $\mA$ is an $m \times n$ matrix,
and $I = \lbrack i_1, i_2, \dots, i_k \rbrack$
and $J = \lbrack j_1, j_2, \dots, j_l \rbrack$,
then $\mA(I, J)$ denotes the $k \times l$ matrix
\begin{equation*}
    \mA(I, J)
    =
    \begin{bmatrix}
        \mA(i_1, j_1) & \mA(i_1, j_2) & \dots & \mA(i_1, j_l) \\
        \mA(i_2, j_1) & \mA(i_2, j_2) & \dots & \mA(i_2, j_l) \\
        \vdots        & \vdots        &       & \vdots        \\
        \mA(i_k, j_1) & \mA(i_k, j_2) & \dots & \mA(i_k, j_l) \\
    \end{bmatrix}
\end{equation*}
We let $\mA(I,:)$ denote the submatrix
$\mA(I, \lbrack 1, 2, \dots, n \rbrack)$
and define $\mA(:, J)$ analogously.
We let $\mA^*$ denote the transpose of $\mA$,
and we say that matrix $\mU$ is \emph{orthonormal}
if its columns are orthonormal, $\mU^* \mU = \mI$.

\subsection{The QR factorization}
The full QR factorization
of a matrix $\mB$ of size $m \times n$
takes the form
\begin{equation}
\label{eq:QR}
\begin{array}{cccc}
    \mB        &=& \mQ        & \mR, \\
    m \times n & & m \times m & m \times n
\end{array}
\end{equation}
where
$\mQ$ is unitary
and
$\mR$ is upper-triangular.

If we further assume that $\mB$ has rank $k$
and that its first $k$ columns are linearly independent,
then it has a rank-$k$ partial QR factorization
given by
\begin{equation*}
\begin{array}{cccc}
    \mB_k      &=& \mQ_k       & \mR_k, \\
    m \times n & & m \times k  & k \times n
\end{array}
\end{equation*}
where $\mQ_k$ has orthonormal columns
and $\mR_k$ is upper-triangular.

\subsection{Randomized compression}
\label{sec:preliminaries_randomized}

Let $\mB$ be an $m \times n$ matrix
that can be accurately approximated by a matrix of rank $k$,
and suppose we seek to determine a matrix $\mQ$
with orthonormal columns (as few as possible) such that
$\Vert \mB - \mQ \mQ^* \mB \Vert$
is small.
In other words, we seek a matrix $\mQ$
whose columns form an approximate orthornomal basis (ON-basis)
for the column space of $\mB$.
This task can efficiently be solved via the following randomized procedure:
\begin{enumerate}
    \item Pick a small integer $p$ representing
        how much ``oversampling'' is done ($p = 10$ is often good),
        and define $r = k + p$.
    \item Form an $n \times r$ matrix $\mG$
        whose entries are independent and identically distributed (i.i.d.)
        normalized Gaussian random numbers.
    \item Form the ``sample matrix'' $\mY = \mB \mG$ of size $m \times r$.
    \item Construct an $m \times r$ matrix $\mQ$ whose columns form an ON basis
        for the columns of $\mY$.
\end{enumerate}

Note that each column of the sample matrix $\mY$
is a random linear combination of the columns of $\mB$.
The probability of the algorithm
producing an accurate result
approaches 1 extremely rapidly as $p$ increases,
and, remarkably, this probability depends only on $p$
(not on $m$ or $n$, or any other properties of $\mB$);
cf.~\cite{halko2011finding}.

\subsection{Functions for computing orthonormal bases}
\label{sec:preliminaries_on_basis}

For matrix $\mB$ of rank at most $k$, we write
\begin{equation*}
\mQ = \col(\mB, k)
\end{equation*}
for a function call that returns a matrix $\mQ$
whose $k$ columns are orthonormal
and span the column space of $\mB$.
In practice, we implement this function using a truncated QR factorization.

For matrix $\mB$ with nullspace of dimension $k$ or greater, we write
\begin{equation*}
\mZ = \nullspace(\mB, k)
\end{equation*}
for a function call that returns matrix $\mZ$
whose $k$ columns are orthonormal
and belong to the nullspace of $\mB$.
In practice, we implement this function
by taking the last $k$ columns of the factor $\mQ$
of the full QR factorization of $\mB^*$.

\begin{remark}
    For a matrix $\mB$ whose first $k$ columns 
    have rank less than $\rank(\mB)$,
    the first $k$ columns of $\mQ$
    produced by an unpivoted QR factorization algorithm
    might not span the column space of $\mB$.
    A similar concern about linear dependence of rows
    exists for the task of computing a basis belonging to the nullspace.
    However, in this work we only apply $\col$ and $\nullspace$
    to random matrices or to products involving random matrices,
    so any subset of $k$ rows or columns of those matrices
    will with probability 1 have rank of $\rank(\mB)$ as long 
    as $\rank(\mB) \leq k$.
    Therefore, $\col$ and $\nullspace$ can safely rely on unpivoted QR factorizations. (Observe that unpivoted QR is numerically stable in this context, and will guarantee that $\|\mB - \mQ\mQ^{*}\mB\|$ is small.)
\end{remark}

\section{Hierarchically block separable matrices} \label{sec:hss}

We start this section with a review of important concepts
relating to HBS matrices,
following the presentation of~\cite{martinsson2019fast}.
While \Cref{sec:tree_structure,sec:hbs_format} review established material
\Cref{sec:telescoping} introduces a new telescoping factorization of HBS matrices
that differs from what is used in prior literature.
We encourage the reader to review section 3.3 carefully,
as the new telescoping factorization is a key idea
leading to the black-box algorithm of \cref{sec:compress}.

\subsection{A tree structure}
\label{sec:tree_structure}

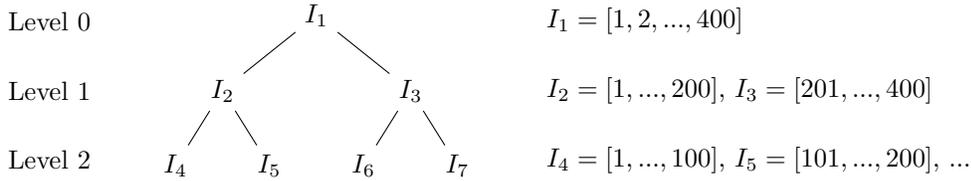
\begin{figure}[tb]
    \centering
    \minipage{0.1\textwidth}
        Level 0
        \vspace{0.5cm}

        Level 1
        \vspace{0.5cm}

        Level 2
    \endminipage
    \hfill
    \minipage{0.35\textwidth}
        \begin{tikzpicture}[level distance=1.0cm,
            level 1/.style={sibling distance=2.5cm},
            level 2/.style={sibling distance=1.25cm}]
            \node{$I_1$}
            child{
                node{$I_2$}
                child{node{$I_4$}}
                child{node{$I_5$}}
            }
            child{
                node{$I_3$}
                child{node{$I_6$}}
                child{node{$I_7$}}
            };
        \end{tikzpicture}
    \endminipage
    \hfill
    \minipage{0.45\textwidth}
        $I_1 = [1, 2, ..., 400]$
        \vspace{0.5cm}

        $I_2 = [1, ..., 200]$,
        $I_3 = [201, ..., 400]$

        \vspace{0.5cm}
        $I_4 = [1, ..., 100]$,
        $I_5 = [101, ..., 200]$, ...
    \endminipage
    \caption{A binary tree structure, where the levels of the tree represent
    successively refined partitions of the index vector $[1, 2, ..., 400]$.
    }
    \label{fig:tree}
\end{figure}

Let $I = \lbrack 1, 2, \dots, N \rbrack$
be a vector of indices corresponding to the rows and columns of an $N \times N$ matrix.
We define a tree $\tree$, in which
each node $\tau$ is associated with a contiguous subset of the indices $I_\tau$.
To the root node of the tree, we assign the full set of indices $I$.
The two children of the root node are given index vectors
$\lbrack 1, \dots, \lceil N/2 \rceil \rbrack$
and
$\lbrack \lceil N/2 \rceil + 1, \dots, N \rbrack$.
We continue evenly splitting the indices of each node to form two child nodes
until we reach a level of the tree
in which the size of each node is below some given threshold $m$.
We refer to a node with children as a parent node,
and a node with no children as a leaf node.
The depth of a node is defined as its distance from the root node,
and level $\ell$ of the tree is defined as the set of nodes with depth $\ell$,
so that level 0 consists of only the root node,
level 1 consists of the two children of the root node, and so on.
The levels of the tree represent successively finer partitions of $I$.
The depth of the tree is defined as the maximum node depth,
denoted by $L \approx \log_2 \left( N / m \right)$.
For simplicity, we only consider fully populated binary trees.
An example tree structure is depicted in \cref{fig:tree}.

\subsection{The HBS matrix format}
\label{sec:hbs_format}

Let $\tree$ be a tree defined on index vector
$I = \lbrack 1, 2, \dots, N \rbrack$.
Matrix $\mA$ of size $N \times N$
is said to be hierarchically block separable
with block rank $k$ with respect to $\tree$
if the following conditions are satisfied.

(1) \emph{Assumption on the ranks of off-diagonal blocks of the finest level:}
For every pair of leaf nodes $\tau$ and $\tau'$,
we define
\begin{equation*}
    \mA_{\tau, \tau'} = \mA(I_\tau, I_{\tau'})
\end{equation*}
and require that if $\tau$ and $\tau'$ are distinct,
then $\mA_{\tau, \tau'}$ must have rank of at most $k$.
Additionally, for each leaf node $\tau$
there must exist basis matrices $\mU_\tau$ and $\mV_\tau$
such that for every leaf node $\tau' \neq \tau$
we have
\begin{equation}
\label{eq:hss_defn_leaf}
\begin{array}{ccccc}
    \mA_{\tau, \tau'} &=& \mU_\tau   & \tilde{\mA}_{\tau,\tau'} & \mV_{\tau'}^*. \\
    m \times m        & & m \times k & k \times k               & k \times m
\end{array}
\end{equation}

(2) \emph{Assumption on the ranks of off-diagonal blocks of levels $L-1, L-2, \dots 1$:}
The following conditions must hold for each level $\ell = L-1, L-2, \dots, 1$.
For every pair of distinct nodes $\tau$ and $\tau'$ on level $\ell$
with children $(\alpha, \beta)$ and $(\alpha', \beta')$, respectively,
we define
\begin{equation*}
    \mA_{\tau, \tau'}
    =
    \begin{bmatrix}
        \tilde{\mA}_{\alpha, \alpha'} & \tilde{\mA}_{\alpha, \beta'} \\
        \tilde{\mA}_{ \beta, \alpha'} & \tilde{\mA}_{ \beta, \beta'} \\
    \end{bmatrix}
\end{equation*}
and require that every such matrix have rank of at most $k$.
Additionally, for each node $\tau$ on level $\ell$
there must exist basis matrices $\mU_\tau$ and $\mV_\tau$
such that for every node $\tau' \neq \tau$ on level $\ell$,
we have
\begin{equation*}
\begin{array}{ccccc}
    \mA_{\tau, \tau'} &=& \mU_\tau    &\tilde{\mA}_{\tau,\tau'}   & \mV_{\tau'}^*. \\
    2k \times 2k      & & 2k \times k & k \times k                & k \times 2k
\end{array}
\end{equation*}

\begin{figure}[tb]
    \centering
    \def\n{8cm}
    \begin{tikzpicture}[level distance=\n/10,
        level 1/.style={sibling distance=\n/2},
        level 2/.style={sibling distance=\n/4},
        level 3/.style={sibling distance=\n/8}]
        \begin{scope}[local bounding box=scope1, line width=1pt]

            \draw (0, 0) rectangle (\n, \n);

            \draw[fill=\lrcolor] (\n*2/4, \n*2/4) rectangle (\n*4/4, \n*4/4) node[midway] {$\mA(I_2, I_3)$};
            \draw[fill=\lrcolor] (\n*0/4, \n*0/4) rectangle (\n*2/4, \n*2/4) node[midway] {$\mA(I_3, I_2)$};

            \draw[fill=\lrcolor] (\n*2/8, \n*6/8) rectangle (\n*4/8, \n*8/8) node[midway] {$\mA(I_4, I_5)$};
            \draw[fill=\lrcolor] (\n*0/8, \n*4/8) rectangle (\n*2/8, \n*6/8) node[midway] {$\mA(I_5, I_4)$};
            \draw[fill=\lrcolor] (\n*6/8, \n*2/8) rectangle (\n*8/8, \n*4/8) node[midway] {$\mA(I_6, I_7)$};
            \draw[fill=\lrcolor] (\n*4/8, \n*0/8) rectangle (\n*6/8, \n*2/8) node[midway] {$\mA(I_7, I_6)$};

            \draw[fill=\lrcolor] (\n*1/8, \n*7/8) rectangle (\n*2/8, \n*8/8) node[midway] {$\mA_{ 8,  9}$};
            \draw[fill=\lrcolor] (\n*3/8, \n*5/8) rectangle (\n*4/8, \n*6/8) node[midway] {$\mA_{10, 11}$};
            \draw[fill=\lrcolor] (\n*5/8, \n*3/8) rectangle (\n*6/8, \n*4/8) node[midway] {$\mA_{12, 13}$};
            \draw[fill=\lrcolor] (\n*7/8, \n*1/8) rectangle (\n*8/8, \n*2/8) node[midway] {$\mA_{14, 15}$};

            \draw[fill=\lrcolor] (\n*0/8, \n*6/8) rectangle (\n*1/8, \n*7/8) node[midway] {$\mA_{ 9,  8}$};
            \draw[fill=\lrcolor] (\n*2/8, \n*4/8) rectangle (\n*3/8, \n*5/8) node[midway] {$\mA_{11, 10}$};
            \draw[fill=\lrcolor] (\n*4/8, \n*2/8) rectangle (\n*5/8, \n*3/8) node[midway] {$\mA_{13, 12}$};
            \draw[fill=\lrcolor] (\n*6/8, \n*0/8) rectangle (\n*7/8, \n*1/8) node[midway] {$\mA_{15, 14}$};

            \draw[fill=\frcolor] (\n*0/8, \n*7/8) rectangle (\n*1/8, \n*8/8) node[midway] {$\mA_{ 8,  8}$};
            \draw[fill=\frcolor] (\n*1/8, \n*6/8) rectangle (\n*2/8, \n*7/8) node[midway] {$\mA_{ 9,  9}$};
            \draw[fill=\frcolor] (\n*2/8, \n*5/8) rectangle (\n*3/8, \n*6/8) node[midway] {$\mA_{10, 10}$};
            \draw[fill=\frcolor] (\n*3/8, \n*4/8) rectangle (\n*4/8, \n*5/8) node[midway] {$\mA_{11, 11}$};
            \draw[fill=\frcolor] (\n*4/8, \n*3/8) rectangle (\n*5/8, \n*4/8) node[midway] {$\mA_{12, 12}$};
            \draw[fill=\frcolor] (\n*5/8, \n*2/8) rectangle (\n*6/8, \n*3/8) node[midway] {$\mA_{13, 13}$};
            \draw[fill=\frcolor] (\n*6/8, \n*1/8) rectangle (\n*7/8, \n*2/8) node[midway] {$\mA_{14, 14}$};
            \draw[fill=\frcolor] (\n*7/8, \n*0/8) rectangle (\n*8/8, \n*1/8) node[midway] {$\mA_{15, 15}$};
        \end{scope}
        \begin{scope}[shift={($(scope1.west)+(-\n*7/20, 0)$)}, grow=right]
            \node{$I_1$}
            child{
                node{$I_3$}
                child{
                    node{$I_7$}
                    child{node{$I_{15}$}}
                    child{node{$I_{14}$}}
                }
                child{
                    node{$I_6$}
                    child{node{$I_{13}$}}
                    child{node{$I_{12}$}}
                }
            }
            child{
                node{$I_2$}
                child{
                    node{$I_5$}
                    child{node{$I_{11}$}}
                    child{node{$I_{10}$}}
                }
                child{
                    node{$I_4$}
                    child{node{$I_{9}$}}
                    child{node{$I_{8}$}}
                }
            };
        \end{scope}
        \begin{scope}[shift={($(scope1.north)+(0, \n*7/20)$)}, grow=down]
            \node{$I_1$}
            child{
                node{$I_2$}
                child{
                    node{$I_4$}
                    child{node{$I_{8}$}}
                    child{node{$I_{9}$}}
                }
                child{
                    node{$I_5$}
                    child{node{$I_{10}$}}
                    child{node{$I_{11}$}}
                }
            }
            child{
                node{$I_3$}
                child{
                    node{$I_6$}
                    child{node{$I_{12}$}}
                    child{node{$I_{13}$}}
                }
                child{
                    node{$I_7$}
                    child{node{$I_{14}$}}
                    child{node{$I_{15}$}}
                }
            };
        \end{scope}
    \end{tikzpicture}
    \caption{
        Tessellation of an HBS matrix with depth 3.
        Low-rank blocks are shown in blue,
        and blocks that are not necessarily low-rank are shown in red.
        The blue blocks are stored in factored form, using ``nested'' basis matrices as described in section
        \ref{sec:hbs_format}.
    }
    \label{fig:matrix_hodlr_lvl3}
\end{figure}
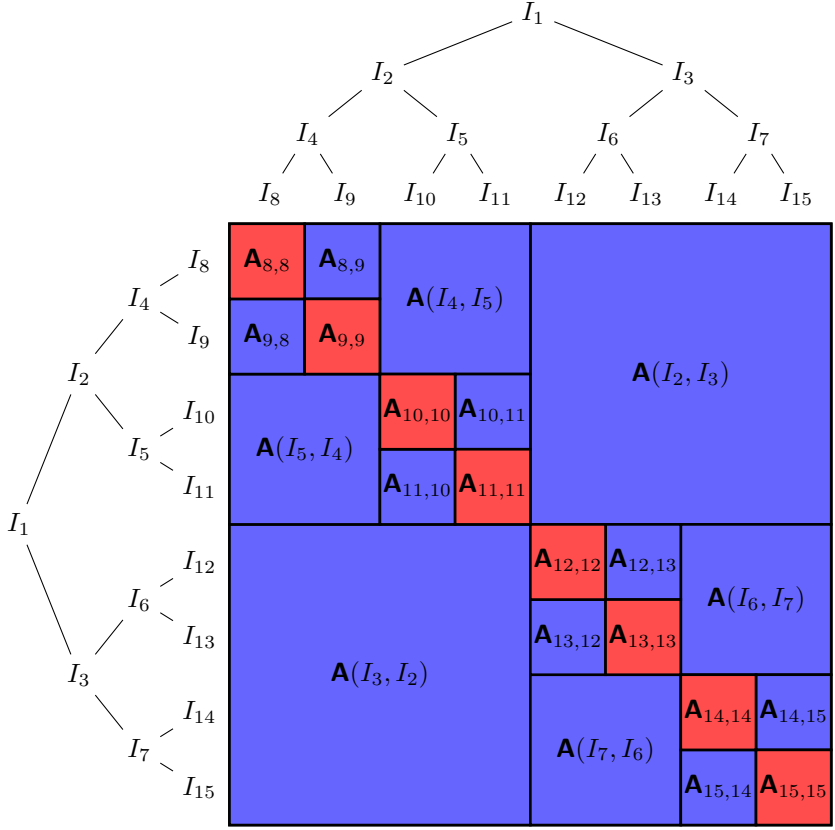

Notably, no assumptions are made on the ranks of the on-diagonal blocks
of level $L$, and those blocks may have full rank.
An example tessellation of an HBS matrix showing compressible and incompressible blocks
is given in \cref{fig:matrix_hodlr_lvl3}.

\subsection{Telescoping factorizations}
\label{sec:telescoping}

We define the following block-diagonal basis matrices.
\begin{equation*}
\begin{aligned}
    \mU^{(\ell)} &= \diag(\mU_\tau : \text{$\tau$ is a node on level $\ell$}), \qquad \ell = 1, 2, \dots, L \\
    \mV^{(\ell)} &= \diag(\mV_\tau : \text{$\tau$ is a node on level $\ell$}), \qquad \ell = 1, 2, \dots, L \\
\end{aligned}
\end{equation*}
Then we obtain a factorization of level $L$ of the form
\begin{equation}
    \label{eq:telescoping_leaf}
    \mA
    =
    \mU^{(L)} \tilde{\mA}^{(L)} (\mV^{(L)})^* + \mD^{(L)},
\end{equation}
where
\begin{align}
\label{eq:newdef1}
    \tilde{\mA}^{(L)} &=  (\mU^{(L)})^* \mA \mV^{(L)} \\
\label{eq:newdef2}
    \mD^{(L)}         &= \mA - \mU^{(L)} \tilde{\mA}^{(L)} (\mV^{(L)})^*.
\end{align}
\Cref{eq:telescoping_leaf} can be viewed as a decomposition of $\mA$
into a term that ``fits'' into the low-rank approximation
using basis matrices $\mU^{(L)}$ and $\mV^{(L)}$
and a discrepancy term $\mD^{(L)}$ that does not.

For successively coarser levels $\ell = L-1, L-2, \dots, 1$,
we similarly have
\begin{equation}
    \label{eq:telescoping_parent}
    \tilde{\mA}^{(\ell+1)}
    =
    \mU^{(\ell)} \tilde{\mA}^{(\ell)} (\mV^{(\ell)})^* + \mD^{(\ell)},
\end{equation}
where
\begin{equation*}
\begin{aligned}
    \tilde{\mA}^{(\ell)} &=  (\mU^{(\ell)})^* \tilde{\mA}^{(\ell+1)} \mV^{(\ell)} \\
    \mD^{(\ell)}         &= \tilde{\mA}^{(\ell+1)} - \mU^{(\ell)} \tilde{\mA}^{(\ell)} (\mV^{(\ell)})^*.
\end{aligned}
\end{equation*}
For the root level, we define
\begin{equation*}
    \mD^{(0)} = \tilde{\mA}^{(\ell+1)}.
\end{equation*}
\Cref{eq:telescoping_leaf,eq:telescoping_parent} define a telescoping factorization of $\mA$.
For example, a factorization with $L=2$ takes the form
\begin{equation*}
    \begin{array}{ccccccccccc}
        \mA &=& \mU^{(2)} & \big(\mU^{(1)} & \mD^{(0)} & (\mV^{(1)})^* &+& \mD^{(1)}\big) & (\mV^{(2)})^* &+& \mD^{(2)}. \\
        \begin{tikzpicture}[scale=1.0, line width=1pt]
            \def\n{13mm}
            \draw (0, 0) rectangle (\n, \n);

            \draw[fill=\lrcolor] (\n*1/4, \n*3/4) rectangle (\n*2/4, \n*4/4);
            \draw[fill=\lrcolor] (\n*2/4, \n*3/4) rectangle (\n*3/4, \n*4/4);
            \draw[fill=\lrcolor] (\n*3/4, \n*3/4) rectangle (\n*4/4, \n*4/4);
            
            \draw[fill=\lrcolor] (\n*0/4, \n*2/4) rectangle (\n*1/4, \n*3/4);
            \draw[fill=\lrcolor] (\n*2/4, \n*2/4) rectangle (\n*3/4, \n*3/4);
            \draw[fill=\lrcolor] (\n*3/4, \n*2/4) rectangle (\n*4/4, \n*3/4);
            
            \draw[fill=\lrcolor] (\n*0/4, \n*1/4) rectangle (\n*1/4, \n*2/4);
            \draw[fill=\lrcolor] (\n*1/4, \n*1/4) rectangle (\n*2/4, \n*2/4);
            \draw[fill=\lrcolor] (\n*3/4, \n*1/4) rectangle (\n*4/4, \n*2/4);
            
            \draw[fill=\lrcolor] (\n*0/4, \n*0/4) rectangle (\n*1/4, \n*1/4);
            \draw[fill=\lrcolor] (\n*1/4, \n*0/4) rectangle (\n*2/4, \n*1/4);
            \draw[fill=\lrcolor] (\n*2/4, \n*0/4) rectangle (\n*3/4, \n*1/4);

            \draw[fill=\frcolor] (\n*0/4, \n*3/4) rectangle (\n*1/4, \n*4/4);
            \draw[fill=\frcolor] (\n*1/4, \n*2/4) rectangle (\n*2/4, \n*3/4);
            \draw[fill=\frcolor] (\n*2/4, \n*1/4) rectangle (\n*3/4, \n*2/4);
            \draw[fill=\frcolor] (\n*3/4, \n*0/4) rectangle (\n*4/4, \n*1/4);
        \end{tikzpicture}
        & &
        \begin{tikzpicture}[scale=1.0, line width=1pt]
            \def\n{13mm}
            \draw (0, 0) rectangle (\n/2, \n);

            \draw[fill=\lrcolor] (\n*0/8, \n*3/4) rectangle (\n*1/8, \n*4/4);
            \draw[fill=\lrcolor] (\n*1/8, \n*2/4) rectangle (\n*2/8, \n*3/4);
            \draw[fill=\lrcolor] (\n*2/8, \n*1/4) rectangle (\n*3/8, \n*2/4);
            \draw[fill=\lrcolor] (\n*3/8, \n*0/4) rectangle (\n*4/8, \n*1/4);
        \end{tikzpicture}
        &
        \begin{tikzpicture}[scale=1.0, line width=1pt]
            \def\n{13mm}
            \fill[white] (0, 0) rectangle (\n/4, \n);
            \draw (0, \n/2) rectangle (\n/4, \n);

            \draw[fill=\lrcolor] (\n*0/8, \n*3/4) rectangle (\n*1/8, \n*4/4);
            \draw[fill=\lrcolor] (\n*1/8, \n*2/4) rectangle (\n*2/8, \n*3/4);
        \end{tikzpicture}
        &
        \begin{tikzpicture}[scale=1.0, line width=1pt]
            \def\n{13mm}
            \fill[white] (0, 0) rectangle (\n/4, \n);
            \draw[fill=\frcolor] (\n*0/4, \n*3/4) rectangle (\n*1/4, \n*4/4);
        \end{tikzpicture}
        &
        \begin{tikzpicture}[scale=1.0, line width=1pt]
            \def\n{13mm}
            \fill[white] (0, 0) rectangle (\n/2, \n);
            \draw (0, \n*3/4) rectangle (\n/2, \n);

            \draw[fill=\lrcolor] (\n*0/4, \n*7/8) rectangle (\n*1/4, \n*8/8);
            \draw[fill=\lrcolor] (\n*1/4, \n*6/8) rectangle (\n*2/4, \n*7/8);
        \end{tikzpicture}
        & &
        \begin{tikzpicture}[scale=1.0, line width=1pt]
            \def\n{13mm}
            \fill[white] (0, 0) rectangle (\n/2, \n);
            \draw[fill=\frcolor] (\n*0/4, \n*3/4) rectangle (\n*1/4, \n*4/4);
            \draw[fill=\frcolor] (\n*1/4, \n*2/4) rectangle (\n*2/4, \n*3/4);
            \draw (0, \n/2) rectangle (\n/2, \n);
        \end{tikzpicture}
        &
        \begin{tikzpicture}[scale=1.0, line width=1pt]
            \def\n{13mm}
            \fill[white] (0, 0) rectangle (\n, \n);
            \draw (0, \n/2) rectangle (\n, \n);

            \draw[fill=\lrcolor] (\n*0/4, \n*7/8) rectangle (\n*1/4, \n*8/8);
            \draw[fill=\lrcolor] (\n*1/4, \n*6/8) rectangle (\n*2/4, \n*7/8);
            \draw[fill=\lrcolor] (\n*2/4, \n*5/8) rectangle (\n*3/4, \n*6/8);
            \draw[fill=\lrcolor] (\n*3/4, \n*4/8) rectangle (\n*4/4, \n*5/8);
        \end{tikzpicture}
        & &
        \begin{tikzpicture}[scale=1.0, line width=1pt]
            \def\n{13mm}
            \draw (0, 0) rectangle (\n, \n);

            \draw[fill=\frcolor] (\n*0/4, \n*3/4) rectangle (\n*1/4, \n*4/4);
            \draw[fill=\frcolor] (\n*1/4, \n*2/4) rectangle (\n*2/4, \n*3/4);
            \draw[fill=\frcolor] (\n*2/4, \n*1/4) rectangle (\n*3/4, \n*2/4);
            \draw[fill=\frcolor] (\n*3/4, \n*0/4) rectangle (\n*4/4, \n*1/4);
        \end{tikzpicture}
    \end{array}
\end{equation*}

It follows from \cref{eq:hss_defn_leaf}
that matrices $\mD^{(\ell)}$, $\ell = 0, 1, \dots L$, are also block-diagonal.
Matrix $\mD^{(\ell)}$ can be described in terms of its on-diagonal blocks as
\begin{equation*}
    \mD^{(\ell)} = \diag(\mD_\tau : \text{$\tau$ is a node on level $\ell$}),
\end{equation*}
where
\begin{equation}
    \label{eq:define_d_tau}
    \mD_\tau = \mA_{\tau, \tau} - \mU_\tau \mU_\tau^* \mA_{\tau, \tau} \mV_\tau \mV_\tau^*.
\end{equation}

\cref{alg:bottom_up_multiply} describes the process
of efficiently applying the telescoping factorization to a vector.

\begin{algorithm}
\caption{Apply a compressed HBS matrix to a vector: $\uu = \mA \qq$.}
\label{alg:bottom_up_multiply}
\begin{algorithmic}
    \State{\underline{Upward pass}}
    \For{level $\ell = L, L-1, \dots, 1$}
        \For{node $\tau$ in level $L$}
            \If{$\tau$ is a leaf node}
                \State{$\hat{\qq}_\tau = \mV_\tau^* \qq(I_\tau)$}
            \Else
                \State{Let $\alpha$ and $\beta$ be the children of $\tau$.}
                \State{$\hat{\qq}_\tau
                    =
                    \mV_\tau^*
                    \begin{bmatrix}
                        \hat{\qq}_\alpha \\
                        \hat{\qq}_\beta
                    \end{bmatrix}
                $}
            \EndIf
        \EndFor
    \EndFor
    \State

    \State{\underline{Downward pass}}
    \For{levels $\ell = 0, 1, \dots, L$}
        \If{$\tau$ is the root node}
            \State{Let $\alpha$ and $\beta$ be the children of $\tau$.}
            \State{$
                \begin{bmatrix}
                    \hat{\uu}_\alpha \\
                    \hat{\uu}_\beta
                \end{bmatrix}
                =
                \mD_\tau
                \begin{bmatrix}
                    \hat{\qq}_\alpha \\
                    \hat{\qq}_\beta
                \end{bmatrix}
            $}
        \ElsIf{$\tau$ is a parent node}
            \State{Let $\alpha$ and $\beta$ be the children of $\tau$.}
            \State{$
                \begin{bmatrix}
                    \hat{\uu}_\alpha \\
                    \hat{\uu}_\beta
                \end{bmatrix}
                =
                \mU_\tau \hat{\uu}_\tau
                +
                \mD_\tau
                \begin{bmatrix}
                    \hat{\qq}_\alpha \\
                    \hat{\qq}_\beta
                \end{bmatrix}
            $}
        \Else
            \State{$\uu(I_\tau) = \mU_\tau \hat{\uu}_\tau + \mD_\tau \qq(I_\tau)$}
        \EndIf
    \EndFor
    \State
\end{algorithmic}
\end{algorithm}

\begin{remark}
\label{rem:factorization_novelty}
Let us stress that our definition of the HBS structure
is slightly different from the more common framework of, e.g.,
\cite{martinsson2011fast,2019_martinsson_fast_direct_solvers}.
Specifically, most authors define the block diagonal matrix
$\mtx{D}^{(L)}$ as the matrix that holds the diagonal subblocks
of the original matrix $\mtx{A}$, which then leaves $\tilde{\mtx{A}}^{(L)}$
with zero blocks on the diagonal.
In contrast, our definition (\ref{eq:newdef2}) has $\mtx{D}^{(L)}$
holding the ``remainder'' of the diagonal blocks after the component
that can be spanned by the basis matrices on level $L$ has been
peeled off in (\ref{eq:newdef1}).
This new definition of $\mtx{D}^{(L)}$ is essential to our
technique for avoiding the need to explicitly form entries of
the original matrix.
\end{remark}

\section{An algorithm for compressing HBS matrices} \label{sec:compress}

This section describes the new compression algorithm for HBS matrices that is
the main contribution of the manuscript.
Let $\mA$ be an $N \times N$ HBS matrix
with block rank $k$ and leaf node size $m$,
and let $r = k + p$, where $p$ represents a small amount of oversampling
($p = 5$ or $p = 10$ are often sufficient).
Let $\mOmega$ and $\mPsi$ be $N \times s$ Gaussian test matrices,
where $s \geq \max(r + m, 3r)$,
and define sample matrices
$\mY = \mA \mOmega$ and $\mZ = \mA^* \mPsi$.
Our objective is to use the information contained in the test and sample matrices
to construct a telescoping factorization of $\mA$,
as defined in \cref{sec:telescoping}.

We begin this section by describing the process
of finding the level-$L$ basis matrices
$\mU^{(L)}$ and $\mV^{(L)}$
and the level-$L$ discrepancy matrix
$\mD^{(L)}$.
Next, we describe how to proceed to coarser levels of the tree.
Finally, we analyze the asymptotic complexity
of the compression algorithm.

\subsection{Computing basis matrices \texorpdfstring{$\mU^{(L)}, \mV^{(L)}$}{U, V}
at the finest level}
\label{sec:compress_basis}

We start by describing the process of computing
the basis matrix $\mU^{(L)}$ of the finest level,
which involves finding for each $\tau$ on level $L$
a basis matrix $\mU_\tau$ that spans the range of
$\mA(I_\tau, I_{\tau'})$ for every node $\tau' \neq \tau$ on level $L$.
We will compute $\mU_\tau$ by applying the randomized algorithm
described in \cref{sec:preliminaries_randomized}
to $\mA(I_\tau, I_\tau^c)$,
where $I_\tau^c = I \setminus I_\tau$
is the set of indices that are not in $I_\tau$.
Importantly, the procedure does not require
the ability to apply $\mA(I_\tau, I_\tau^c)$
to random vectors;
rather we compute the randomized samples using 
only information contained in $\mOmega$ and $\mY$.

For the purpose of illustration,
suppose for now that $\mA$ is an HBS matrix with depth 2,
and we want to find the basis matrix $\mU_4$
associated with node 4.
We let
$\mOmega_4 = \mOmega(I_4, :)$
and
$\mY_4 = \mY(I_4, :)$
denote the blocks of $\mOmega$ and $\mY$
of size $m \times s$ associated with node 4.
Then the randomized sample of $\mA$
takes the following form.
\begin{equation*}
    \begin{array}{cccc}
        \mY &=& \mA & \mOmega \\
        \begin{tikzpicture}[scale=1.0, line width=1pt]
            \def\n{20mm}
            \draw[fill=\randcolor] (\n*0/8, \n*3/4) rectangle (\n*3/8, \n*4/4) node[midway] {$\mY_4$};
            \draw[fill=\randcolor, dashed] (\n*0/8, \n*2/4) rectangle (\n*3/8, \n*3/4) node[midway] {$\mY_5$};
            \draw[fill=\randcolor, dashed] (\n*0/8, \n*1/4) rectangle (\n*3/8, \n*2/4) node[midway] {$\mY_6$};
            \draw[fill=\randcolor, dashed] (\n*0/8, \n*0/4) rectangle (\n*3/8, \n*1/4) node[midway] {$\mY_7$};
        \end{tikzpicture}
        & &
        \begin{tikzpicture}[scale=1.0, line width=1pt]
            \def\n{20mm}
            \draw[fill=\lrcolor] (\n*1/4, \n*3/4) rectangle (\n*2/4, \n*4/4);
            \draw[fill=\lrcolor] (\n*2/4, \n*3/4) rectangle (\n*3/4, \n*4/4);
            \draw[fill=\lrcolor] (\n*3/4, \n*3/4) rectangle (\n*4/4, \n*4/4);

            \draw[fill=\lrcolor, dashed] (\n*0/4, \n*2/4) rectangle (\n*1/4, \n*3/4);
            \draw[fill=\lrcolor, dashed] (\n*2/4, \n*2/4) rectangle (\n*3/4, \n*3/4);
            \draw[fill=\lrcolor, dashed] (\n*3/4, \n*2/4) rectangle (\n*4/4, \n*3/4);

            \draw[fill=\lrcolor, dashed] (\n*0/4, \n*1/4) rectangle (\n*1/4, \n*2/4);
            \draw[fill=\lrcolor, dashed] (\n*1/4, \n*1/4) rectangle (\n*2/4, \n*2/4);
            \draw[fill=\lrcolor, dashed] (\n*3/4, \n*1/4) rectangle (\n*4/4, \n*2/4);

            \draw[fill=\lrcolor, dashed] (\n*0/4, \n*0/4) rectangle (\n*1/4, \n*1/4);
            \draw[fill=\lrcolor, dashed] (\n*1/4, \n*0/4) rectangle (\n*2/4, \n*1/4);
            \draw[fill=\lrcolor, dashed] (\n*2/4, \n*0/4) rectangle (\n*3/4, \n*1/4);

            \draw[fill=\frcolor] (\n*0/4, \n*3/4) rectangle (\n*1/4, \n*4/4);
            \draw[fill=\frcolor, dashed] (\n*1/4, \n*2/4) rectangle (\n*2/4, \n*3/4);
            \draw[fill=\frcolor, dashed] (\n*2/4, \n*1/4) rectangle (\n*3/4, \n*2/4);
            \draw[fill=\frcolor, dashed] (\n*3/4, \n*0/4) rectangle (\n*4/4, \n*1/4);
        \end{tikzpicture}
        &
        \begin{tikzpicture}[scale=1.0, line width=1pt]
            \def\n{20mm}
            \draw[fill=\randcolor] (\n*0/8, \n*3/4) rectangle (\n*3/8, \n*4/4) node[midway] {$\mOmega_4$};
            \draw[fill=\randcolor] (\n*0/8, \n*2/4) rectangle (\n*3/8, \n*3/4) node[midway] {$\mOmega_5$};
            \draw[fill=\randcolor] (\n*0/8, \n*1/4) rectangle (\n*3/8, \n*2/4) node[midway] {$\mOmega_6$};
            \draw[fill=\randcolor] (\n*0/8, \n*0/4) rectangle (\n*3/8, \n*1/4) node[midway] {$\mOmega_7$};
        \end{tikzpicture}
        \\
    \end{array}
\end{equation*}

Since $\mOmega_4$ is of size $m \times s$,
it has a nullspace of dimension at least $s - m \geq r$.
Then we can find a set of $r$ orthonormal vectors
that belong to its nullspace,
$\mP_4 = \nullspace(\mOmega_4, r)$,
so that $\mOmega_4 \mP_4 = \mzero$.
Then we have
\begin{equation*}
    \begin{array}{cccc}
        \mY \mP_4 &=& \mA & \mOmega \mP_4. \\
        \begin{tikzpicture}[scale=1.0, line width=1pt]
            \def\n{20mm}
            \draw[fill=\randcolor] (\n*0/8, \n*3/4) rectangle (\n*1/8, \n*4/4);
            \draw[fill=\randcolor, dashed] (\n*0/8, \n*2/4) rectangle (\n*1/8, \n*3/4);
            \draw[fill=\randcolor, dashed] (\n*0/8, \n*1/4) rectangle (\n*1/8, \n*2/4);
            \draw[fill=\randcolor, dashed] (\n*0/8, \n*0/4) rectangle (\n*1/8, \n*1/4);
        \end{tikzpicture}
        & &
        \begin{tikzpicture}[scale=1.0, line width=1pt]
            \def\n{20mm}
            \draw[fill=\lrcolor] (\n*1/4, \n*3/4) rectangle (\n*2/4, \n*4/4);
            \draw[fill=\lrcolor] (\n*2/4, \n*3/4) rectangle (\n*3/4, \n*4/4);
            \draw[fill=\lrcolor] (\n*3/4, \n*3/4) rectangle (\n*4/4, \n*4/4);
            
            \draw[fill=\lrcolor, dashed] (\n*0/4, \n*2/4) rectangle (\n*1/4, \n*3/4);
            \draw[fill=\lrcolor, dashed] (\n*2/4, \n*2/4) rectangle (\n*3/4, \n*3/4);
            \draw[fill=\lrcolor, dashed] (\n*3/4, \n*2/4) rectangle (\n*4/4, \n*3/4);
            
            \draw[fill=\lrcolor, dashed] (\n*0/4, \n*1/4) rectangle (\n*1/4, \n*2/4);
            \draw[fill=\lrcolor, dashed] (\n*1/4, \n*1/4) rectangle (\n*2/4, \n*2/4);
            \draw[fill=\lrcolor, dashed] (\n*3/4, \n*1/4) rectangle (\n*4/4, \n*2/4);
            
            \draw[fill=\lrcolor, dashed] (\n*0/4, \n*0/4) rectangle (\n*1/4, \n*1/4);
            \draw[fill=\lrcolor, dashed] (\n*1/4, \n*0/4) rectangle (\n*2/4, \n*1/4);
            \draw[fill=\lrcolor, dashed] (\n*2/4, \n*0/4) rectangle (\n*3/4, \n*1/4);

            \draw[fill=\frcolor] (\n*0/4, \n*3/4) rectangle (\n*1/4, \n*4/4);
            \draw[fill=\frcolor, dashed] (\n*1/4, \n*2/4) rectangle (\n*2/4, \n*3/4);
            \draw[fill=\frcolor, dashed] (\n*2/4, \n*1/4) rectangle (\n*3/4, \n*2/4);
            \draw[fill=\frcolor, dashed] (\n*3/4, \n*0/4) rectangle (\n*4/4, \n*1/4);
        \end{tikzpicture}
        &
        \begin{tikzpicture}[scale=1.0, line width=1pt]
            \def\n{20mm}
            \draw[fill=white] (\n*0/8, \n*0/4) rectangle (\n*1/8, \n*4/4);
            \draw[fill=\randcolor] (\n*0/8, \n*2/4) rectangle (\n*1/8, \n*3/4);
            \draw[fill=\randcolor] (\n*0/8, \n*1/4) rectangle (\n*1/8, \n*2/4);
            \draw[fill=\randcolor] (\n*0/8, \n*0/4) rectangle (\n*1/8, \n*1/4);
        \end{tikzpicture}
        \\
    \end{array}
\end{equation*}
Notably, submatrix $(\mOmega \mP_4)(I_4,:)$
is filled with zeros,
and submatrix $(\mOmega \mP_4)(I_4^c,:)$
is filled with values that turn out to also have
a standard Gaussian distribution,
cf.~\cref{remark:whygaussian}.
The product $\mA \mOmega \mP_4$
can be viewed as a randomized sample of $\mA$,
excluding contributions from columns $\mA(:, I_4)$,
so the rows of $\mA \mOmega \mP_4$
indexed by $I_4$ contain a randomized sample
of $\mA(I_4, I_4^c)$.
As suggested by the diagram above,
we can obtain that sample inexpensively
by simply multiplying $\mY_4 \mP_4$.
Then we orthonormalize the sample
to find the basis matrix,
$\mU_4 = \col(\mY_4 \mP_4, r)$.

The generalization to an arbitrary leaf node $\tau$
is straightforward.
We define
$$
\mOmega_\tau = \mOmega(I_\tau, :)
\qquad\mbox{and}\qquad
\mY_\tau     = \mY(I_\tau, :).
$$
Then we compute the nullspace
\begin{equation*}
    \mP_\tau = \nullspace(\mOmega_\tau, r)
\end{equation*}
so that $\mY_\tau \mP_\tau$
contains a randomized sample
of $\mA(I_\tau, I_\tau^c)$.
Finally, we orthonormalize the sample
to find the basis matrix $\mU_\tau$,
\begin{equation}
\label{eq:u_tau_formula}
    \mU_\tau = \col(\mY_\tau \mP_\tau, r).
\end{equation}

A similar process using $\mPsi$ and $\mZ$
yields basis matrices $\mV_\tau$,
\begin{equation}
\label{eq:v_tau_formula}
    \begin{aligned}
        \mQ_\tau &= \nullspace(\mPsi_\tau, r) \\
        \mV_\tau &= \col(\mZ_\tau \mQ_\tau, r).
    \end{aligned}
\end{equation}

\begin{remark}
\label{remark:whygaussian}
We claimed that the matrix $(\mOmega \mP_\tau)(I_\tau^c, :) = 
\mOmega(I_\tau^c, :) \mP_\tau$ 
has a Gaussian distribution. This claim follows from (i) the
fact that the distribution of Gaussian matrices is invariant
under unitary transformations, and (ii) that the construction 
of the matrix $\mP_\tau$ is done independently of the entries in 
$\mOmega(I_\tau^c, :)$.
\end{remark}

\subsection{Computing discrepancy matrix \texorpdfstring{$\mD^{(L)}$}{D}
at the finest level}
\label{sec:compress_discrepancy}

Once we have computed $\mU^{(L)}$ and $\mV^{(L)}$,
we next compute $\mD^{(L)}$.
Recall that our definition of $\mD^{(L)}$,
given in~\cref{eq:define_d_tau},
is different from what appears
in the typical telescoping factorization of HBS matrices
(cf.~\cref{rem:factorization_novelty}).
We proceed by rewriting~\cref{eq:define_d_tau} as
\begin{equation}
    \label{eq:d_tau_sum}
    \mD_\tau = (\mI - \mU_\tau \mU_\tau^*) \mA_{\tau, \tau} + \mU_\tau \mU_\tau^* \mA_{\tau, \tau} (\mI - \mV_\tau \mV_\tau^*),
\end{equation}
and deriving formulas for computing
$(\mI - \mU_\tau \mU_\tau^*) \mA_{\tau, \tau}$
and
$\mU_\tau \mU_\tau^* \mA_{\tau, \tau} (\mI - \mV_\tau \mV_\tau^*)$
separately.

For $(\mI - \mU_\tau \mU_\tau^*) \mA_{\tau, \tau}$,
we first express $\mY_\tau$ as a blocked matrix product
\begin{equation*}
    \mY_\tau = \sum_\text{$\tau'$ in level $\ell$} \mA_{\tau, \tau'} \mOmega_{\tau'}.
\end{equation*}
Multiplying $(\mI - \mU_\tau \mU_\tau^*)$ and applying \cref{eq:hss_defn_leaf} gives
\begin{equation*}
    (\mI - \mU_\tau \mU_\tau^*) \mY_\tau
    =
    (\mI - \mU_\tau \mU_\tau^*) \mA_{\tau, \tau} \mOmega_\tau.
\end{equation*}
Solving a least-squares problem with $\mOmega_\tau$ gives
\begin{equation}
    \label{eq:d_tau_part_1}
    (\mI - \mU_\tau \mU_\tau^*) \mA_{\tau, \tau}
    =
    (\mI - \mU_\tau \mU_\tau^*) \mY_\tau \mOmega_\tau^\dagger.
\end{equation}

A similar derivation yields
\begin{equation}
    \label{eq:d_tau_part_2}
    \mA_{\tau, \tau} (\mI - \mV_\tau \mV_\tau^*)
    =
    ( (\mI - \mV_\tau \mV_\tau^*) \mZ_\tau \mPsi_\tau^\dagger)^*.
\end{equation}
Substituting~\cref{eq:d_tau_part_1,eq:d_tau_part_2}
into~\cref{eq:d_tau_sum}
gives the formula
\begin{equation}
\label{eq:d_tau_formula}
    \mD_\tau
    =
    (\mI - \mU_\tau \mU_\tau^*) \mY_\tau \mOmega_\tau^\dagger
    +
    \mU_\tau \mU_\tau^* ( (\mI - \mV_\tau \mV_\tau^*) \mZ_\tau \mPsi_\tau^\dagger)^*.
\end{equation}

\begin{remark}
    The computation of $\mD_\tau$
    involves solving least squares problems
    with Gaussian matrices $\mOmega_\tau$ and $\mPsi_\tau$
    of size $m \times s$ for $\tau$ belonging to level $L$
    or size $2r \times s$ for $\tau$ belonging to a coarser level.
    Gaussian matrices with nearly square shapes
    have non-negligible probabilities
    of being ill-conditioned,
    but the probabilities quickly become negligible
    even for slightly rectangular matrices~\cite{chen2005condition,edelman2005tails}.
    Such concerns can be alleviated
    by choosing $s$ to be sufficiently large.
    For the numerical experiments in \cref{sec:experiments},
    we simply use $s = r + m = 3r$.
\end{remark}

\subsection{Compressing levels \texorpdfstring{$\ell = L-1, L-2, \dots, 0$}{above level L}}
\label{sec:compress_coarse_levels}

After compressing level $L$,
we proceed to the next coarser level $L-1$.
That is, we seek
$\mU^{(L-1)}, \mV^{(L-1)},$ and $\mD^{(L-1)}$
that satisfy~\cref{eq:telescoping_parent}.
We will first obtain randomized samples of $\tilde{\mA}^{(L)}$,
and then using the same procedure as for level $L$
we will find $\mU^{(L-1)}, \mV^{(L-1)},$ and $\mD^{(L-1)}$.

To compute randomized samples of $\tilde{\mA}^{(L)}$,
we multiply~\cref{eq:telescoping_leaf} with $\mOmega$ to obtain
\begin{equation*}
    \mY
    =
    \mA \mOmega
    =
    (\mU^{(L)} \tilde{\mA}^{(L)} (\mV^{(L)})^* + \mD^{(L)}) \mOmega,
\end{equation*}
and rearrange to obtain
\begin{equation}
\label{eq:newtestmatrix}
    \underbrace{(\mU^{(L)})^* (\mY - \mD^{(L)} \mOmega)}_\text{sample matrix}
    = \tilde{\mA}^{(L)}
    \underbrace{(\mV^{(L)})^* \mOmega}_\text{test matrix}.
\end{equation}
Then the columns of $(\mU^{(L)})^* (\mY - \mD^{(L)} \mOmega)$
contain $s$ samples of
$\tilde{\mA}^{(L)}$
taken with test matrix $(\mV^{(L)})^* \mOmega$.
Then for node $\tau$ on level $L-1$
with children $\alpha$ and $\beta$,
we define
\begin{equation}
\label{eq:omega_tau_and_y_tau}
    \mOmega_\tau =
        \begin{bmatrix}
            \mV_\alpha^* \mOmega_\alpha \\
            \mV_\beta^* \mOmega_\beta
        \end{bmatrix}
    \qquad \text{and} \qquad
    \mY_\tau =
        \begin{bmatrix}
        \mU_\alpha^* (\mY_\alpha - \mD_\alpha \mOmega_\alpha) \\
        \mU_\beta^*  (\mY_\beta  - \mD_\beta  \mOmega_\beta)
        \end{bmatrix},
\end{equation}
to be the corresponding blocks of the new test and sample matrices.
We define $\mPsi_\tau$ and $\mZ_\tau$ analogously.
Once we have $\mOmega_\tau, \mPsi_\tau, \mY_\tau, \mZ_\tau$
we compute $\mU_\tau, \mV_\tau,$ and $\mD_\tau$
exactly as before using \cref{eq:u_tau_formula,eq:v_tau_formula,eq:d_tau_formula}.

This process is applied to successively coarser levels of the tree
until the root node is reached.
For the root node $\tau$, we have
$\mY_\tau = \mD^{(0)} \mOmega_\tau$,
so we simply solve
$\mD^{(0)} = \mY_\tau \mOmega_\tau^\dagger$.
The full compression procedure is summarized in \cref{alg:bottom_up_compress}.

\begin{remark}

The sampling at the finest level that was described in \cref{sec:compress_basis} is directly supported by the theory for the randomized SVD \cite{halko2011finding,martinsson2020randomized} since the test matrices are provably Gaussian matrices drawn independently of the matrices that are approximated, cf.~Remark \ref{remark:whygaussian}.
The sampling at the higher levels is harder to analyze.
Considering equation \cref{eq:newtestmatrix}, we see that the test matrix for sampling $\tilde{\mA}^{(L)}$ is $(\mV^{(L)})^* \mOmega$,
and for level $\ell$, the test matrix for sampling
$\tilde{\mA}^{(\ell)}$
is
$(\mV^{(\ell)})^* (\mV^{(\ell+1)})^* \cdots (\mV^{(L)})^* \mOmega$.
Since for $\ell = L-1, ..., 1$,
$\mV^{(\ell)}$ depends not only on $\mA$ and $\mPsi$, but also on $\mOmega$,
the argument in \cref{remark:whygaussian} does not immediately apply. 
Extensive numerical experiments (including those reported in \cref{sec:experiments}) indicate that these test matrices work just as well as ``true'' Gaussian ones, but we have not yet been able to substantiate this claim through analysis.
\end{remark}

\begin{algorithm}
\caption{Compressing an HBS matrix}
\label{alg:bottom_up_compress}
\begin{algorithmic}
    \State{\underline{Compute randomized samples of $\mA$ and $\mA^*$.}}
    \State{Form Gaussian random test matrices $\mOmega$ and $\mPsi$
        of size $N \times s$.}
    \State{Multiply $\mY = \mA \mOmega$ and $\mZ = \mA^* \mPsi$.}
    \Statex

    \State{\underline{Compress level by level from finest to coarsest.}}
    \For{level $\ell = L, L-1, \dots, 0$}
        \For{node $\tau$ in level $\ell$}
            \State{\underline{Obtain test and sample matrices associated with $\tau$.}}
            \If{$\tau$ is a leaf node}
                \State{
                    $\begin{aligned}
                    \mOmega_\tau &= \mOmega(I_\tau, :),
                        \quad
                        &\mPsi_\tau &= \mPsi(I_\tau, :) \\
                    \mY_\tau &= \mY(I_\tau, :),
                        \quad
                        &\mZ_\tau &= \mZ(I_\tau, :)
                    \end{aligned}$
                }
            \Else
                \State{Let $\alpha$ and $\beta$ denote the children of $\tau$.}
                \State{
                    $\begin{aligned}
                    \mOmega_\tau &=
                        \begin{bmatrix}
                        \mV_\alpha^* \mOmega_\alpha \\
                        \mV_\beta^* \mOmega_\beta
                        \end{bmatrix},
                    \quad
                    &\mPsi_\tau &=
                        \begin{bmatrix}
                        \mU_\alpha^* \mPsi_\alpha \\
                        \mU_\beta^* \mPsi_\beta
                        \end{bmatrix} \\
                    \mY_\tau &=
                        \begin{bmatrix}
                        \mU_\alpha^* (\mY_\alpha - \mD_\alpha \mOmega_\alpha) \\
                        \mU_\beta^*  (\mY_\beta  - \mD_\beta  \mOmega_\beta)
                        \end{bmatrix},
                    \quad
                    &\mZ_\tau &=
                        \begin{bmatrix}
                        \mV_\alpha^* (\mZ_\alpha - \mD_\alpha^* \mPsi_\alpha) \\
                        \mV_\beta^*  (\mZ_\beta  - \mD_\beta^*  \mPsi_\beta)
                        \end{bmatrix}
                    \end{aligned}$
                }
            \EndIf
            \Statex
            \State{\underline{Compute blocks of the factorization associated with $\tau$.}}
            \If{$\ell > 0$}
                \State{$\begin{aligned}
                    \mP_\tau &= \nullspace(\mOmega_\tau, r),
                    \quad
                    &\mQ_\tau &= \nullspace(\mPsi_\tau, r) \\
                    \mU_\tau &= \col(\mY_\tau \mP_\tau, r),
                    \quad
                    &\mV_\tau &= \col(\mZ_\tau \mQ_\tau, r)
                    \end{aligned}$}

                \State{$\mD_\tau
                    =
                    (\mI - \mU_\tau \mU_\tau^*) \mY_\tau \mOmega_\tau^\dagger
                    +
                    \mU_\tau \mU_\tau^* \left( (\mI - \mV_\tau \mV_\tau^*) \mZ_\tau \mPsi_\tau^\dagger \right)^*
                    $}
            \Else
                \State{$\mD_\tau = \mY_\tau \mOmega_\tau^\dagger$}
                \Comment{Here $\tau$ is the root node.}
            \EndIf
        \EndFor
    \EndFor
\end{algorithmic}
\end{algorithm}

\subsection{Asymptotic complexity}

Recall that $r = k + p$,
where $k$ is the block rank of $\mA$
and $p$ represents the amount of oversampling.
We assume for simplicity that the leaf node size is $m = 2r$
and $\mOmega$ and $\mPsi$ both have $s = 3r$ columns.
\cref{alg:bottom_up_compress} requires $s$ matrix-vector products
of $\mA$ and $\mA^*$,
and an additional $\bigO(r^3)$ operations
for each node in the tree,
of which there are approximately $2 N / m$.
Therefore, the total compression time is
\begin{equation*}
    \tcompress = 6 r N \times \trand + 6 r \times \tmult + \bigO(r^2 N) \times \tflop,
\end{equation*}
where $\trand$ denotes the time to sample a value from the standard Gaussian distribution,
$\tmult$ denotes the time to apply $\mA$ or $\mA^*$ to a vector,
and $\tflop$ denotes the time to carry out a floating point arithmetic operation.

\begin{remark}[Comparison of information efficiency]
    If we view each randomized sample
    as carrying $N$ values worth of information about $\mA$,
    and still assume $m = 2r$ and $s = 3r$,
    we find that the compression algorithm requires a total of $6rN$ values
    to reconstruct the matrix.
    The algorithm of~\cite{martinsson2011fast}
    requires only $r$ samples of $\mA$ and $\mA^*$,
    but it also requires access to $\sim mN$ matrix entries
    that form the block-diagonal part of $\mA$
    as well as $\sim rN$ elements that appear in interpolative decompositions,
    for a total of $5rN$ values worth of information.
    Therefore, the algorithm in the present work
    requires only slightly more information to recover $\mA$,
    while having the advantage of being truly black box.
\end{remark}

\section{Numerical experiments} \label{sec:experiments}
In this section, we present a selection of numerical results.
We report the following quantities
for a number of test problems:
(1) the time to compress the operator,
(2) the time to apply the compressed representation to a vector,
(3) the relative accuracy of the compressed representation,
and (4) the storage requirements of the compressed representation
measured as the number of floating point
values per degree of freedom.
The algorithms for compressing matrices and applying compressed representations
are written in Python,
and the black-box multiplication routines are written in MATLAB.
The experiments were carried out on a workstation
with an Intel Core i9-10900K processor with 10 cores
and 128GB of memory.

We report two measurements of compression time:
(1) the time spent executing
\cref{alg:bottom_up_compress},
and (2) the execution time of
\cref{alg:bottom_up_compress},
excluding time spent within the black-box multiplication routine.
In the context of black-box compression,
it is assumed that the black-box multiplication routine is provided,
and any steps taken to initialize the matvec routine
are considered external to the compression algorithm.
Therefore, the reported timings exclusively measure time spent executing
\cref{alg:bottom_up_compress},
and do not include any time initializing the matvec routines.
For each of the test problems we consider in this section,
the time spent initializing the matvec routines is relatively short,
taking less than one tenth of the time spent for compression.

We measure the accuracy of the compressed matrices using the relative error
\[
\frac{\Vert \tilde{\mA} - \mA \Vert}{\Vert \mA \Vert},
\]
where $\mtx{A}$ is the matrix defined by the matvec in each example.
The operator norms are estimated by taking 20 steps of power iteration.
Specifically, to estimate the norm of
$\mB \in \{\mA, \tilde{\mA} - \mA\}$,
we initialize a random vector $\vv_0$,
compute $\vv_i = \mB^* \mB (\vv_{i-1} / \Vert \vv_{i-1} \Vert)$,
$i = 1, ..., 20$,
and take the estimate
$\Vert \mB \Vert \approx \sqrt{\Vert \vv_{20} \Vert}$.
We also report the maximum leaf node size $m$
and the number $r$ of random vectors per test matrix,
which are inputs to the compression algorithm.
We select the value of $r$
so that the resulting approximation
achieves reasonably high accuracy
and then set $m = 2r$ and $s = 3r$.

\subsection{Boundary integral equation}
\label{sec:bie}
\begin{figure}
\begin{center}
\includegraphics[height=30mm]{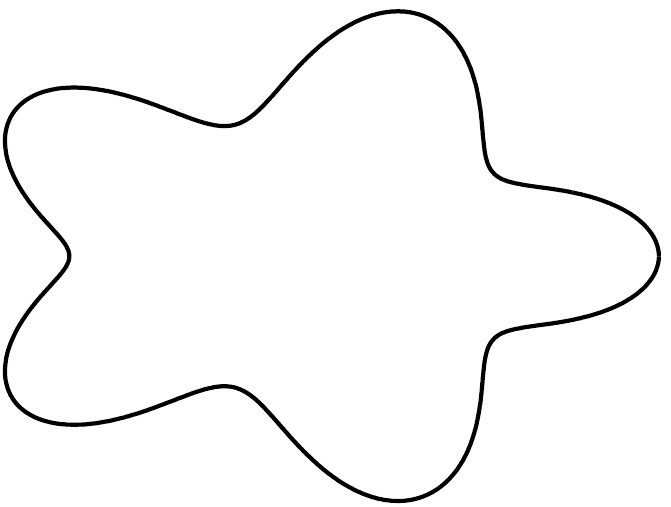}
\end{center}
\caption{Contour $\Gamma$ on which the BIE (\ref{eq:BIE}) is defined.}
\label{fig:BIE}
\end{figure}

We consider a matrix arising from the discretization of the
Boundary Integral Equation (BIE)
\begin{equation}
\label{eq:BIE}
\frac{1}{2}q(\pvct{x}) +
\int_{\Gamma}
\frac{(\pvct{x} - \pvct{y})\cdot\pvct{n}(\pvct{y})}
     {4\pi|\pvct{x} - \pvct{y}|^{2}}\,q(\pvct{y})\,ds(\pvct{y}) = f(\pvct{x}),
\qquad\pvct{x} \in \Gamma,
\end{equation}
where $\Gamma$ is the simple closed contour in the plane shown in Figure \ref{fig:BIE},
and where
$\pvct{n}(\pvct{y})$ is the outwards pointing unit normal of $\Gamma$ at
$\pvct{y}$. The BIE (\ref{eq:BIE}) is a standard integral equation formulation
of the Laplace equation with boundary condition $f$ on the domain interior
to $\Gamma$. The BIE (\ref{eq:BIE}) is discretized using the
Nystr\"om method on $N$ equispaced points on $\Gamma$, with the Trapezoidal
rule as the quadrature (since the kernel in (\ref{eq:BIE}) is smooth,
the Trapezoidal rule has exponential convergence).

The fast matrix-vector multiplication is in this case furnished by the recursive skeletonization (RS) procedure of \cite{2005_martinsson_fastdirect}.
In this case, the ``exact'' matrix $\mtx{A}$ is the RS approximation to the discretized integral operator.
To avoid spurious effects due to the rank structure inherent 
in the RS representation, we set a computational tolerance in the RS approximation to $1e-15$. 
Additionally, we used a different tree structure for the RS compression than the one we used when running the randomized algorithm.

\begin{figure}[tb]
    \begin{subfigure}{0.48\textwidth}
        \includegraphics[width=\textwidth]{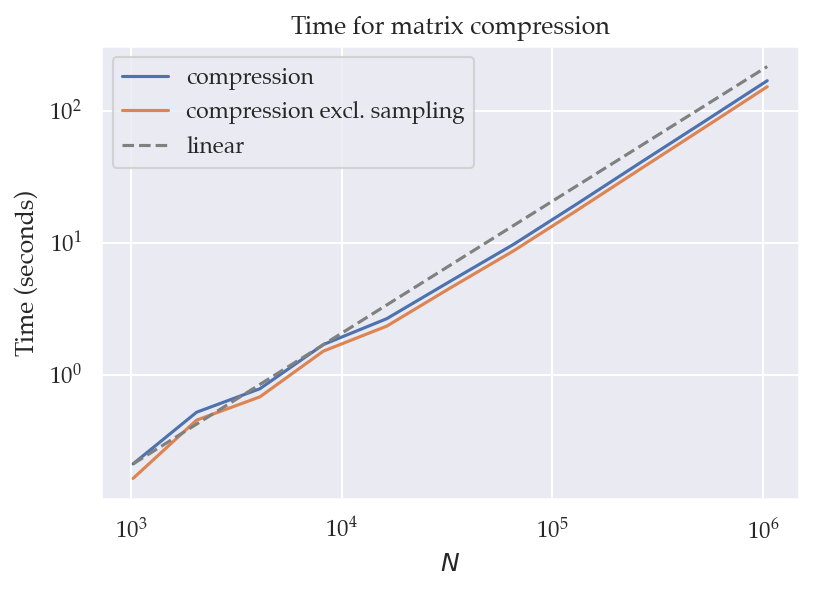}
    \end{subfigure}
    \hspace*{\fill}
    \begin{subfigure}{0.48\textwidth}
        \includegraphics[width=\textwidth]{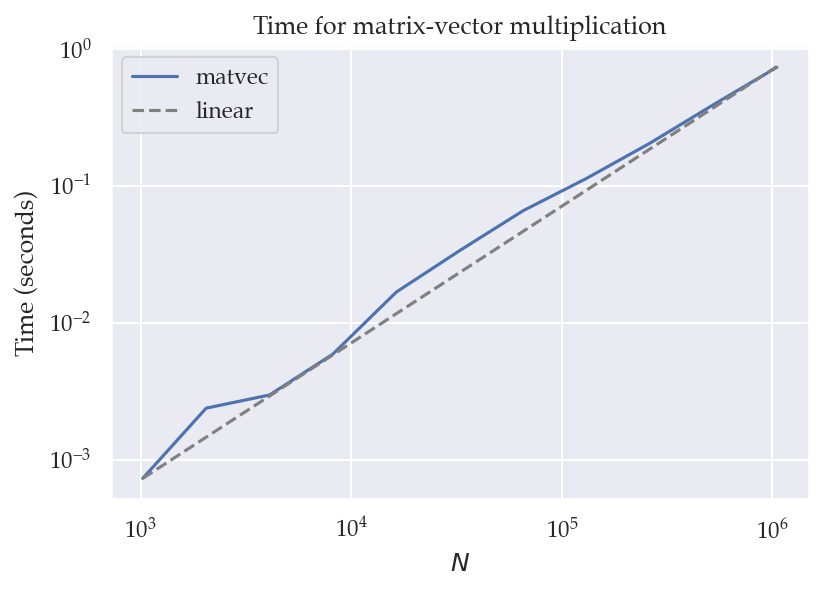}
    \end{subfigure}

    \medskip
    \begin{subfigure}{0.48\textwidth}
        \includegraphics[width=\textwidth]{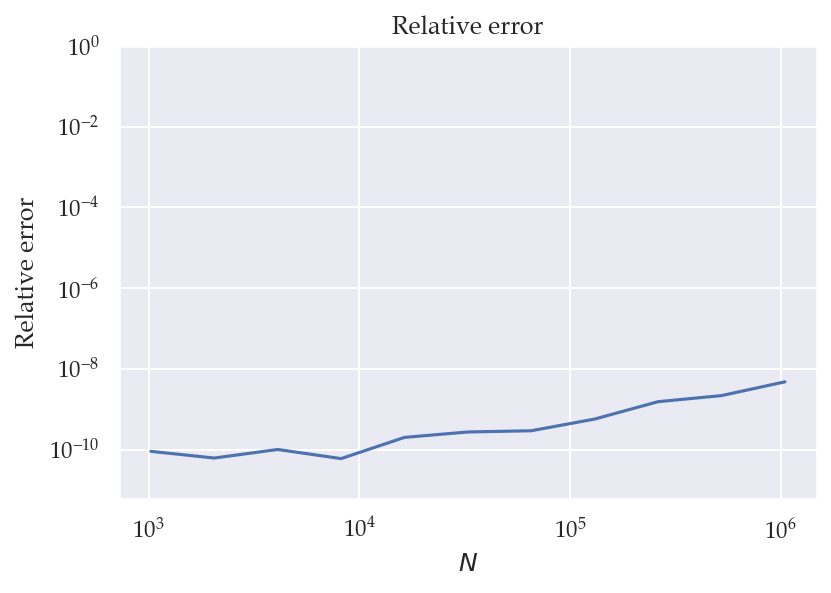}
    \end{subfigure}
    \hspace*{\fill}
    \begin{subfigure}{0.48\textwidth}
        \includegraphics[width=\textwidth]{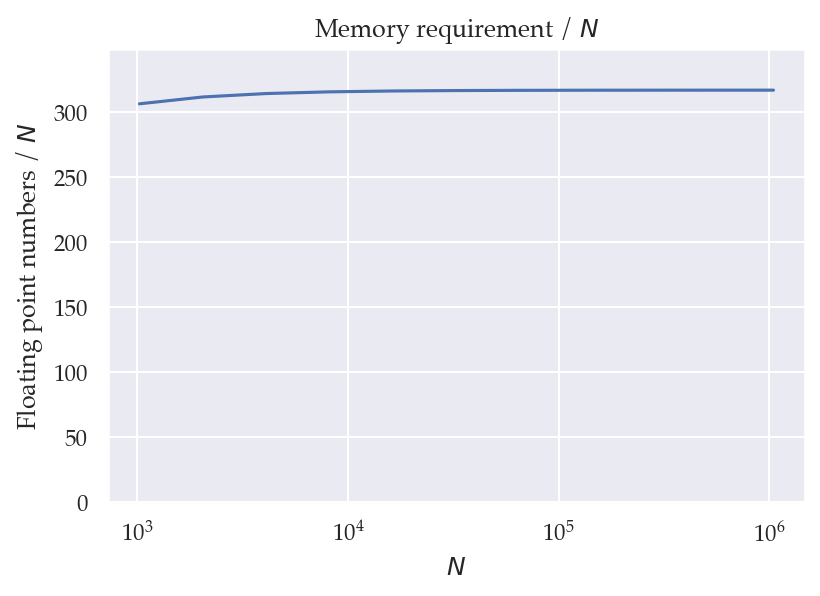}
    \end{subfigure}
    \caption{
        Results from applying the compression algorithm
        to a double layer potential on a simple contour in the plane.
        Here $r = 30$ and $m = 60$.
    }
    \label{fig:dl_results}
\end{figure}

Results are given in \cref{fig:dl_results}.

\begin{remark}
The problem under consideration here is artificial in the sense that there is
no actual need to use more than a couple of hundred points to resolve (\ref{eq:BIE})
numerically to double precision accuracy. It is included merely to illustrate
the asymptotic performance of the proposed method.
\end{remark}

\subsection{Operator multiplication}
\label{sec:operator_mult}

\begin{figure}[tb]
    \begin{subfigure}{0.48\textwidth}
        \includegraphics[width=\textwidth]{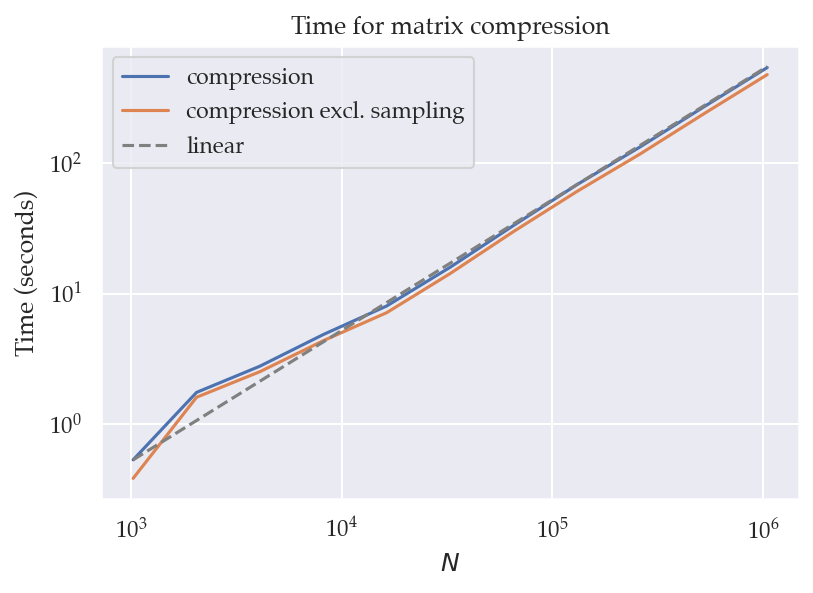}
    \end{subfigure}
    \hspace*{\fill}
    \begin{subfigure}{0.48\textwidth}
        \includegraphics[width=\textwidth]{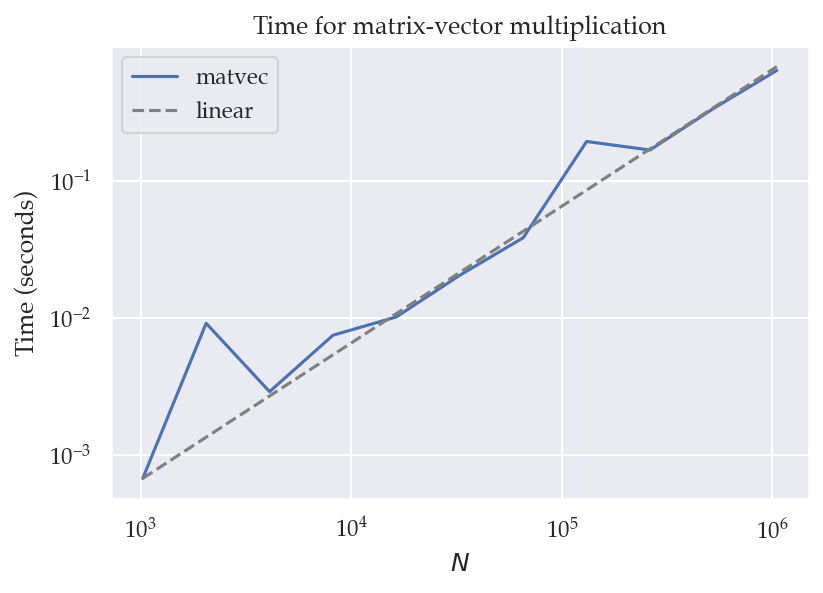}
    \end{subfigure}

    \medskip
    \begin{subfigure}{0.48\textwidth}
        \includegraphics[width=\textwidth]{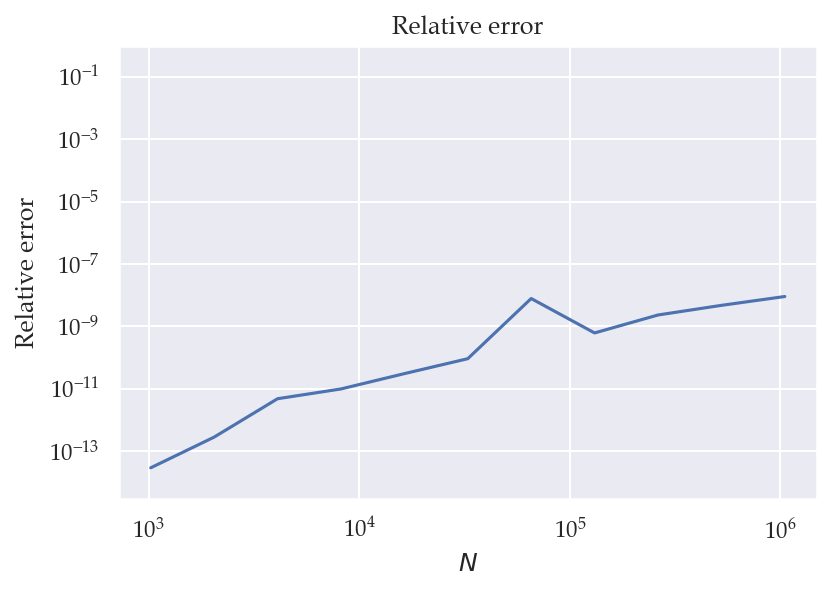}
    \end{subfigure}
    \hspace*{\fill}
    \begin{subfigure}{0.48\textwidth}
        \includegraphics[width=\textwidth]{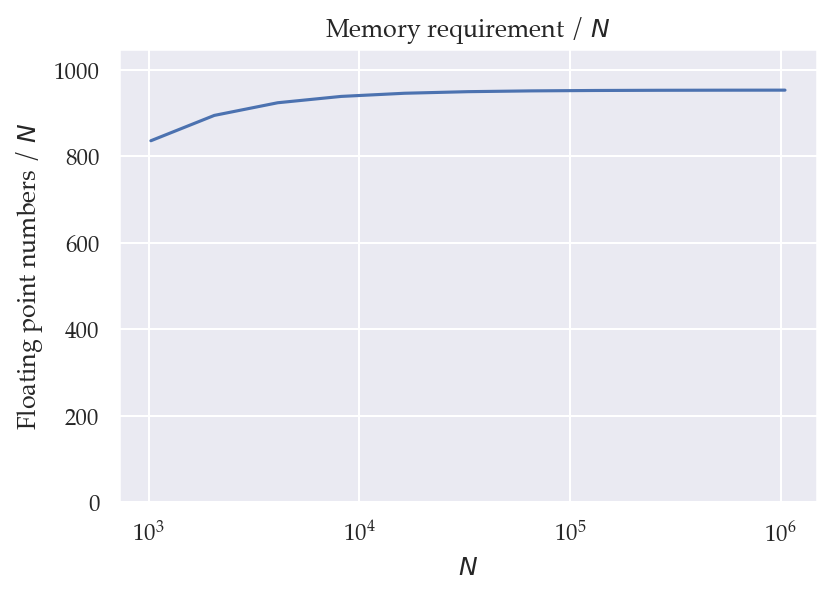}
    \end{subfigure}
    \caption{
        Results from applying the compression algorithm
        to the Neumann-to-Dirichlet operator.
        Here $r = 100$ and $m = 200$.
    }
    \label{fig:dfn_results}
\end{figure}

We next investigate how the proposed technique performs on a matrix matrix
multiplication problem. Specifically, we determine the Neumann-to-Dirichlet
operator $T$ for the contour shown in Figure \ref{fig:BIE} using the well
known formula
\[
T = S\left(\frac{1}{2}I + D^*\right)^{-1},
\]
where $S$ is the single layer operator
$
[Sq](\pvct{x}) = \int_{\Gamma}-\frac{1}{2\pi}\,\log|\pvct{x} - \pvct{y}|\,q(\pvct{y})\,ds(\pvct{y}),
$
and where $D^{*}$ is the adjoint of the double-layer operator
$
[D^{*}q](\pvct{x}) = \int_{\Gamma}\frac{\pvct{n}(\pvct{x}) \cdot (\pvct{x} - \pvct{y})}{2\pi |\pvct{x} - \pvct{y}|^{2}}
\,q(\pvct{y})\,ds(\pvct{y}).
$
The operators $S$ and $D$ are again discretized using a Nystr\"om method
on equispaced points (with sixth order Kapur-Rokhlin \cite{1997_kapur_rokhlin}
corrections to handle the singularity in $S$), resulting in matrices
$\mtx{S}$ and $\mtx{D}$. The matrix $\mtx{S}\bigl(\tfrac{1}{2}\mtx{I} + \mtx{D}^{*}\bigr)^{-1}$
is again applied using the recursive skeletonization procedure of
\cite{2005_martinsson_fastdirect}.
Results are given in \cref{fig:dfn_results}.

\subsection{Frontal matrices in nested dissection}
\label{sec:nesteddiss}
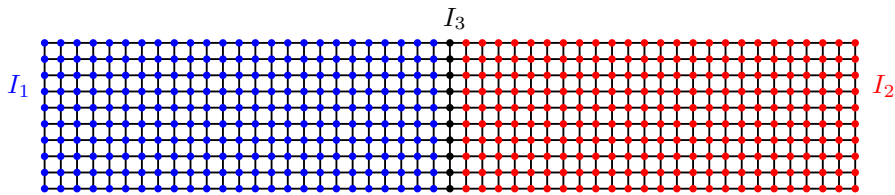
\begin{figure}
\begin{center}
\setlength{\unitlength}{1mm}
\begin{picture}(125,26)
\put(005,00){\includegraphics[width=110mm]{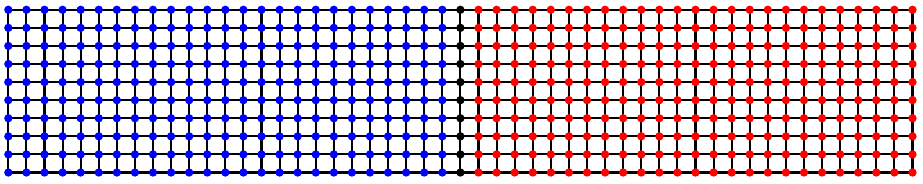}}
\put(001,14){$\color{blue}I_{1}$}
\put(116,14){$\color{red}I_{2}$}
\put(059,23){$\color{black}I_{3}$}
\end{picture}
\end{center}
\caption{An example of the grid in the sparse LU example
described in Section \ref{sec:nesteddiss}.
There are $N \times n$ points in the grid, shown for $N = 10, n = 51$.
}
\label{fig:nesteddiss_geom}
\end{figure}

Our next example is a simple model problem that illustrates the behavior
of the proposed method in the context of sparse direct solvers. The idea
here is to use rank structure to compress the increasingly large Schur
complements that arise in the LU factorization of a sparse matrix arising
from the finite element or finite difference discretization of an
elliptic PDE, cf.~\cite[Ch.~21]{2019_martinsson_fast_direct_solvers}.
As a model problem, we consider
a $51N \times 51N$ matrix
$\mtx{C}$ that
encodes the stiffness matrix for the standard five-point stencil
finite difference approximation to the Poisson equation on a rectangle
using a grid with $N \times 51$ nodes. We partition the grid
into three sets $\{1,2,3\}$, as shown in \cref{fig:nesteddiss_geom},
and then tessellate $\mtx{C}$ accordingly,
\[
\mtx{C} =
\begin{bmatrix}
\mtx{C}_{11} & \mtx{0}      & \mtx{C}_{13} \\
\mtx{0}      & \mtx{C}_{22} & \mtx{C}_{23} \\
\mtx{C}_{31} & \mtx{C}_{32} & \mtx{C}_{33}
\end{bmatrix},
\]
where 
$\mtx{C}_{11}$
and
$\mtx{C}_{22}$
are $25N \times 25N$,
and
$\mtx{C}_{33}$
is $N \times N$.
The matrix we seek to compress is the $N \times N$ Schur complement
\[
\mtx{A} = \mtx{C}_{33}
- \mtx{C}_{31}\mtx{C}_{11}^{-1}\mtx{C}_{13}
- \mtx{C}_{32}\mtx{C}_{22}^{-1}\mtx{C}_{23}.
\]
In our example, we apply $\mtx{A}$ to vector by calling standard
sparse direct solvers for the left and the right subdomains, respectively.

\begin{figure}[tb]
    \begin{subfigure}{0.48\textwidth}
        \includegraphics[width=\textwidth]{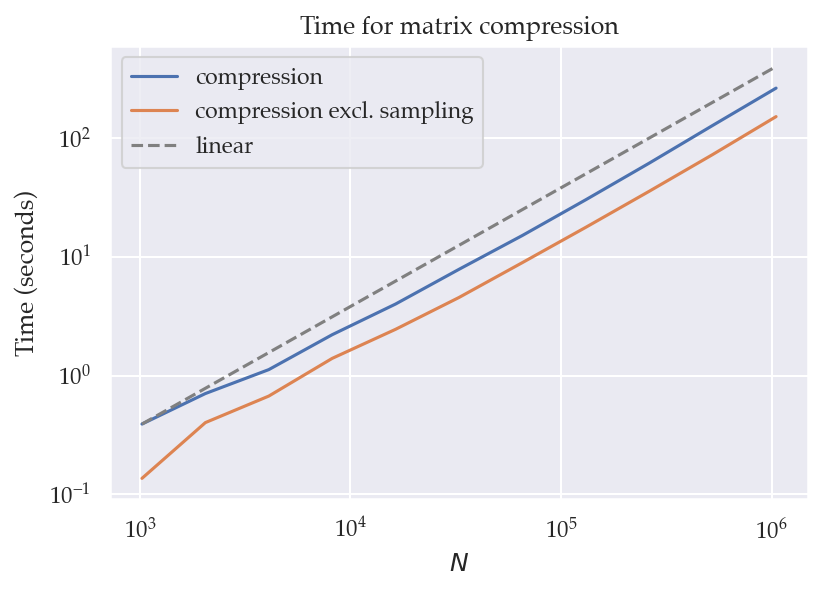}
    \end{subfigure}
    \hspace*{\fill}
    \begin{subfigure}{0.48\textwidth}
        \includegraphics[width=\textwidth]{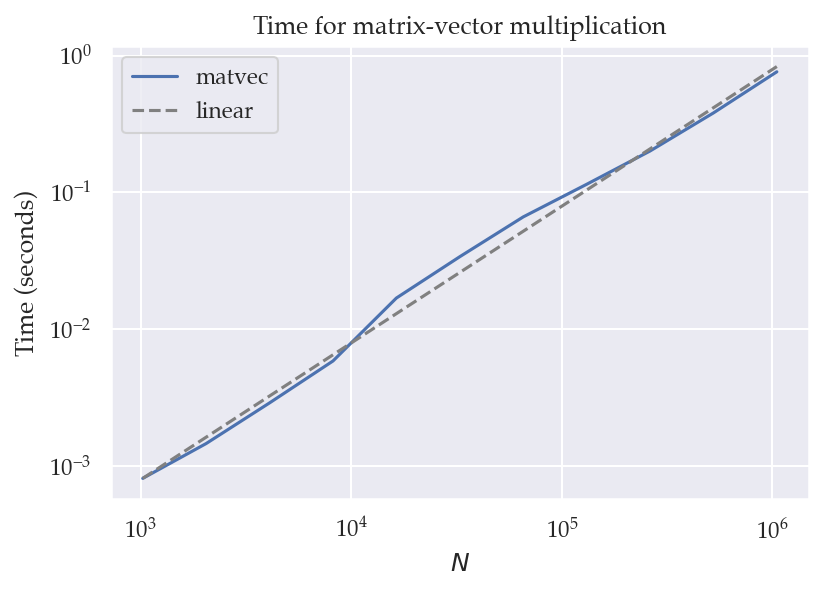}
    \end{subfigure}

    \medskip
    \begin{subfigure}{0.48\textwidth}
        \includegraphics[width=\textwidth]{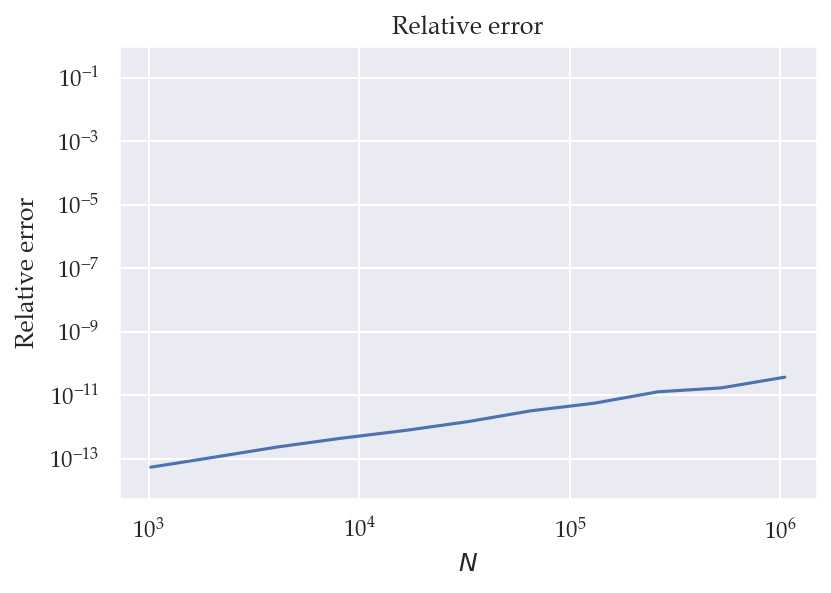}
    \end{subfigure}
    \hspace*{\fill}
    \begin{subfigure}{0.48\textwidth}
        \includegraphics[width=\textwidth]{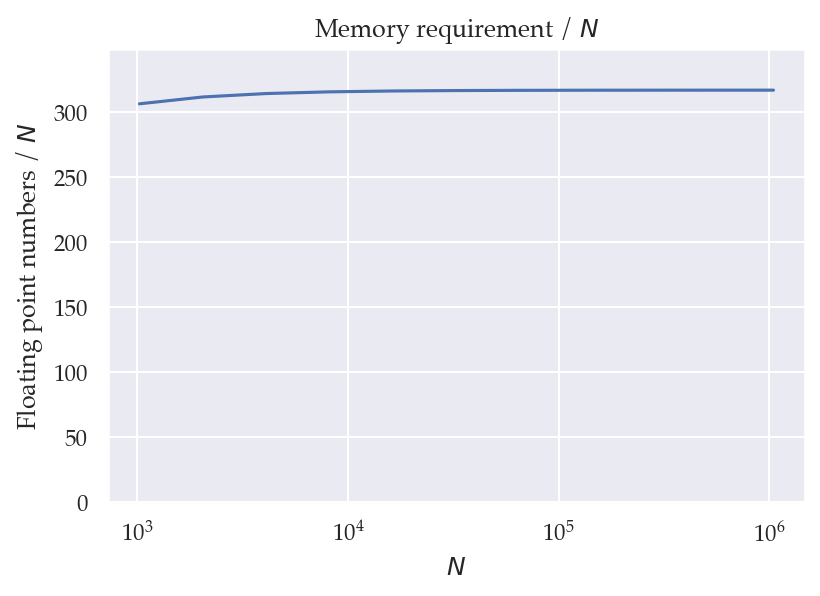}
    \end{subfigure}
    \caption{
        Results from applying the compression algorithm
        to frontal matrices in the nested dissection algorithm.
        Here $r = 30$ and $m = 60$.
    }
    \label{fig:nd_results}
\end{figure}

Results are given in \cref{fig:nd_results}.

\subsection{Summary of observations}
\begin{itemize}
\item The numerical results support our claim of linear complexity for all steps of the algorithms presented.
\item The approximations achieve high accuracy in every case.
In particular, the numerical results indicate that the method
is numerically stable, and that there is minimal aggregation
of errors as hierarchical trees deepen.
\item The asymptotic storage cost,
    reported as number of floating point numbers per degree of freedom,
    appears to be independent of the problem size $N$.
    (It does, of course, depend on the rank parameter $r$ and on the size of the diagonal blocks $m$.)
\end{itemize}

\section{Conclusions} \label{sec:conclusions}
This paper presents an algorithm
for black-box randomized compression
of Hierarchically Block Separable matrices.
To compress an $N \times N$ matrix $\mA$,
the algorithm requires
only $\bigO(k)$ samples of $\mA$ and $\mA^*$,
where $k$ is the block rank of $\mA$.
Numerical experiments demonstrate that the algorithm
is accurate and very computationally efficient,
with compression time scaling linearly in $N$
when the cost of applying $\mA$ and $\mA^*$
to a vector is $\bigO(N)$.

\subsection*{Acknowledgments}
The work reported was supported by the Office of Naval Research (N00014-18-1-2354),
by the National Science Foundation
(DMS-2313434 and DMS-1952735),
and by the Department of Energy ASCR (DE-SC0022251).

\appendix

\bibliographystyle{siamplain}
\bibliography{references}
\end{document}